\documentclass[final]{article}
\usepackage[latin9]{inputenc}
\usepackage{fullpage}
\usepackage{amsmath,amsfonts,amsthm}
\usepackage{url,xspace,calc,ifthen}

\usepackage{listings}
\lstdefinelanguage{pseudo}{basicstyle=\ttfamily,keywordstyle=\bfseries,%
                           keywords={Input,Output,for,end,if,then,else,procedure,function,return},%
                           escapechar=!}
\newcommand{\svd}{\textsc{svd}\xspace} 
\newcommand{\qr}{\textsc{qr}\xspace}   
\newcommand{\lapack}{\textsc{lapack}\xspace} 

\renewcommand{\vector}[1]{\ensuremath{\mathbf{#1}}}
\renewcommand{\matrix}[1]{\ensuremath{\mathbf{#1}}}

  \newcommand{\R}{\ensuremath{\mathbb{R}}}
\renewcommand{\H}{\ensuremath{\mathbb{H}}}

\newcommand{\norm}[1]{\ensuremath{\|#1\|}}    
\newcommand{\transpose}[1]{\ensuremath{#1^T}} 
\newcommand{\hermitian}[1]{\transpose{\conjugate{#1}}} 
\newcommand{\conjugate}[1]{\ensuremath{\overline{#1}}} 
\newcommand{\inverse}[1]{\ensuremath{#1^{-1}}} 

\newtheorem{theorem}{Theorem}
\newtheorem{lemma}{Lemma}

\newcommand{\range}{\to} 

\title{Quaternion Singular Value Decomposition based on Bidiagonalization
       to a Real Matrix using Quaternion Householder Transformations.}

\author{S.~J. Sangwine\footnotemark[2]\ \footnotemark[4]
 \and   N. Le~Bihan\footnotemark[3]} 

\begin{document}
\maketitle
\renewcommand{\thefootnote}{\fnsymbol{footnote}}
\footnotetext[2]{Department of Electronic Systems Engineering, University of Essex,
                 Wivenhoe Park, Colchester CO4 3SQ, United Kingdom. \texttt{s.sangwine@ieee.org}}
\footnotetext[3]{Laboratoire des Images et des Signaux, CNRS
                 UMR 5083, ENSIEG, 961 Rue de la Houille Blanche,
                 Domaine Universitaire, BP 46, 38402 Saint Martin d'Hères,
                 Cedex, France. \texttt{nicolas.le-bihan@lis.inpg.fr}}
\footnotetext[4]{The work presented here was carried out at the Laboratoire des Images et des Signaux,
                 Grenoble. Financial support from the Royal Academy of Engineering of the United Kingdom
                 and the Centre National de la Recherche Scientifique (CNRS) in France, which
                 enabled Dr~Sangwine to work in Grenoble for 7 months is gratefully acknowledged.}


\begin{abstract}
We present a practical and efficient means to compute the singular value decomposition (\svd) of a
quaternion matrix \matrix{A} based on bidiagonalization of \matrix{A} to a \emph{real} bidiagonal
matrix \matrix{B} using quaternionic Householder transformations. Computation of the \svd of
\matrix{B} using an existing subroutine library such as \lapack provides the singular values of
\matrix{A}. The singular vectors of \matrix{A} are obtained trivially from the product of the
Householder transformations and the real singular vectors of \matrix{B}. We show in the paper that
left and right quaternionic Householder transformations are different because of the non-commutative
multiplication of quaternions and we present formulae for computing the Householder vector and
matrix in each case.
\end{abstract}


\section{Introduction}
The singular value decomposition (\svd) of a quaternion matrix was first described theoretically in
1997 by Zhang~\cite[Theorem 7.2]{Zhang:97}. Zhang demonstrated the existence of the quaternion \svd
using an isomorphism between a quaternion matrix and a complex matrix known as the \emph{complex
adjoint}. The details of the complex adjoint matrix need not concern us here, but it is important to
note that this matrix is redundant --- it contains twice the number of real values as the
corresponding quaternion matrix. That is, every quaternion is represented by four complex numbers in
the complex adjoint, rather than the two that are sufficient, the pattern of sign/conjugation and
the position of these complex numbers preserving the properties of the quaternion, so that computing
with the complex adjoint is equivalent to computing with the quaternion matrix. Therefore, computing
the \svd of the complex adjoint matrix is equivalent to computing the \svd of the corresponding
quaternion matrix, but since the \svd of the complex adjoint contains the singular values of the
quaternion matrix twice, alternate singular values of the complex adjoint have to be discarded.
Details on using the complex adjoint matrix to compute the quaternion \svd may be found in
\cite{BihanMars:2004,LeBihan:2001} but this method is effectively superseded by the results in this
paper\footnote{It is also possible to construct a real `adjoint' of a quaternion matrix, but computing
with it suffers from the same problems of redundancy as computing with a complex adjoint.}.

It has been noted by the authors of a quaternion \qr algorithm~\cite{Bunse-Gerstner:1989} that using
complex adjoint matrices to compute quaternion matrix decompositions can lead to loss of accuracy,
since an algorithm applied to a complex adjoint does not necessarily preserve accurately the pattern
of values in the complex matrix. For this reason, and because the computation inherently involves
twice as much work as a direct quaternion computation, we have sought direct quaternion methods to
compute the quaternion \svd.

In this paper, we present a method to compute the \svd of an arbitrary quaternion matrix \matrix{A}
which is both elegant and fast, based on transformation of \matrix{A} to \emph{real} bidiagonal form
using quaternionic Householder transformations. We show in the paper that quaternionic Householder
transformations have two forms, dependent on whether the Householder matrix is to be applied on the
left or right of the vector to be transformed (this is because quaternion multiplication is not
commutative). Using left and right Householder transformations, the bidiagonalization of a
quaternionic matrix may be represented by: $\matrix{H}\matrix{A}\matrix{G} = \matrix{B}$, where
\matrix{H} (\matrix{G}) is the product of the left (right) quaternion Householder matrices required
to transform \matrix{A} to the real bidiagonal matrix \matrix{B}. Since \matrix{H} and \matrix{G}
are unitary, we have also\footnote{In this paper we use the overbar to denote a quaternion conjugate
and a superscript $T$ to denote a transpose. We avoid the standard superscript asterisk in our work,
reserving it to denote a complex conjugate in work with complexified quaternions.}: $\matrix{A} =
\hermitian{\matrix{H}}\matrix{B}\hermitian{\matrix{G}}$. The \svd of the real bidiagonal matrix
\matrix{B} may be computed using existing algorithms and code for real matrices. This is usually
represented by $\matrix{B} = \matrix{U}\matrix{\Sigma}\transpose{\matrix{V}}$, where \matrix{U} and
\matrix{V} are orthogonal real matrices. We show that the \svd of the quaternion matrix \matrix{A}
is then:
\begin{equation}
\matrix{A}=\hermitian{\matrix{H}}\matrix{U}\matrix{\Sigma}\transpose{\matrix{V}}\hermitian{\matrix{G}}
\end{equation}
Thus the singular values of the quaternion matrix \matrix{A} are the same as the singular values of
the bidiagonal real matrix \matrix{B}, and the singular vectors of \matrix{A} are obtained trivially
from the products $\hermitian{\matrix{H}}\matrix{U}$ and
$\transpose{\matrix{V}}\hermitian{\matrix{G}}$.

The most significant aspect of this method is that the \svd of the real matrix can be computed using
an existing subroutine or library. Considerable effort has been expended by others (for example, the
authors of the \lapack library~\cite{LAPACKUG}) to compute the \svd of real (and complex) matrices
efficiently and accurately and the results of this effort can be exploited in computing quaternion
\svd{}s without having to write quaternion versions of the \svd computation itself. The idea of
computing the \svd using reduction to bidiagonal form is well known, but it is much less well known
that the bidiagonal form may be real, even for a complex matrix\footnote{Interested readers may like
to inspect the code of the \lapack routine \textsc{zgebrd} or refer to the relevant section of the
\lapack Users Guide~\cite{LAPACKUG} available online at
\url{http://www.netlib.org/lapack/lug/node53.html}.}.

The sequence of topics in the rest of the paper is as follows. In §\ref{quaternionhouseholder} we
discuss the generalization of the Householder transformation to quaternion vectors; we prove that
the transformation exists, and is unitary for both the left and right cases; and we present an
algorithm for computing the quaternionic Householder vector in each case. In §\ref{existence} we
prove the existence of the real bidiagonalization of an arbitrary quaternion matrix using the
quaternionic Householder transformations. In §\ref{bidiag} we give an algorithm for bidiagonalizing
a quaternion matrix. We then show in §\ref{qsvd} that the singular values of the bidiagonalized
matrix are identical to the singular values of the original quaternion matrix, and we present an
algorithm for computing the quaternion \svd using the quaternion bidiagonalization in §\ref{bidiag}
and any existing algorithm for the \svd of a real matrix.

We mention here for completeness that it is also possible to compute the real bidiagonal matrix
\matrix{B} using a quaternionic form of the Golub-Kahan-Lanczos
bidiagonalization~\cite{GolubKahan:1965,Golub:1996}. This is conceptually simpler than the use of
Householder transformations, since the only change needed to generalize the algorithm to handle
quaternion matrices is the definition of the norm of a quaternion vector, which we give here in
Lemma \ref{vectornorm}. However, the Golub-Kahan-Lanczos algorithm is not attractive because of its
computational complexity. For each non-zero element of the bidiagonal matrix \matrix{B} generated,
the corresponding quaternion vector is generated, but to yield accurate results, these quaternion
vectors must be re-orthogonalized against all the previously generated vectors, and this requires
considerable computation. We have implemented the Golub-Kahan bidiagonalization for quaternion
matrices and verified that it works, but the method presented in this paper (which is also a
generalization of the Householder method first presented in \cite{GolubKahan:1965}) is much faster,
because it does not require re-orthogonalization.

We assume familiarity with quaternions and quaternion matrices. The general concepts of quaternions
are available in several books~\cite{Kuipers:1999, Ward:1997, Kantor:1989}. The paper by Zhang
\cite{Zhang:97} is currently the most convenient and comprehensive summary of material on quaternion
matrices. Some elementary results which we need in the rest of the paper are given in Appendix
\ref{elementary}.


\section{Quaternion Householder transformations}
\label{quaternionhouseholder}
The Householder transformation~\cite{DOI:10.1145/320941.320947} is based on a unitary matrix of the
form $\matrix{H} = \left(\matrix{I}-\vector{u}\transpose{\vector{u}}\right)$, where
$\norm{\vector{u}}=\sqrt{2}$. \vector{u} is called a Householder vector. \matrix{H} is used to
transform a vector, either by multiplication on the left by \matrix{H} for a column vector, or by
multiplication on the right by \transpose{\matrix{H}} for a row vector. Normally the vectors are
parts of a row or column of a matrix that is transformed by a succession of Householder
transformations. The modification usually consists of setting to zero multiple elements of the
vector, while preserving the norm, but the general case is interesting. Given a unit vector
\vector{v} and a row or column vector \vector{a} to be transformed, the transformation effectively
generates a copy of \vector{v} scaled by the norm of \vector{a}, that is:
$\matrix{H}\vector{a}=\norm{\vector{a}}\vector{v}$. Thus if \vector{v} contains one non-zero element
with a modulus of 1, the transformed vector will have one non-zero element in the corresponding
position with a modulus equal to the norm of \vector{a} and sign equal to the sign of the non-zero
element of \vector{v}.

The Householder transformation was generalised to the complex case by
Morrison~\cite{DOI:10.1145/321021.321030}. He introduced a complex scalar with unit modulus,
denoted~$\zeta$. The Householder matrix and vector then become complex: $\matrix{H} =
z\left(\matrix{I}-\vector{u}\hermitian{\vector{u}}\right)$ where $z = 1/\zeta$. Now, if \vector{v}
has non-zero \emph{real} elements, and a unit norm, then the transformed vector will be zero in all
positions except those where \vector{v} was non-zero, and the non-zero elements of the transformed
vector will be real. This property follows from the fact that the transformation generates a copy of
\vector{v} scaled by the norm of \vector{a} --- if the elements of \vector{v} are imaginary or
complex, then so will be the corresponding elements of the transformed vector. Note that by using a
real vector \vector{v} the transformation results in a \emph{real} modified vector. Therefore a
sequence of Householder transformations may be used to bidiagonalize a complex matrix directly to a
real bidiagonal matrix.

In the complex case, $z$ commutes with elements of the matrix
$\matrix{I}-\vector{u}\hermitian{\vector{u}}$ (and with elements of \vector{a}). Therefore, the
same matrix is used whether multiplied on the right or left (but as noted above, if multiplied
on the right it must be transposed).

In the quaternion case, we must distinguish between a \emph{left} Householder transformation (to be
applied to a column vector) and a \emph{right} Householder transformation (to be applied to a row
vector). To make this distinction clear, we use a different symbol in each case, \matrix{H} for a
left Householder matrix, and \matrix{G} for a right Householder matrix. In computing the Householder
matrix, we must be careful to compute and use $z$ and \vector{u} correctly, because they do not
commute with elements of \vector{a}. It is worth emphasising that the terms \emph{left} and
\emph{right} are used in the matrix multiplication sense and not in the quaternion multiplication
sense, although the difference between the two cases arises from the non-commutativity of
quaternion multiplication.

In this section, we generalize the Householder transformation to the case of an arbitrary quaternion
vector \vector{a}, but we are forced to restrict the vector \vector{v} to be real in order to solve
for the Householder vector \vector{u} and the scalar $\zeta$. This is not an important restriction,
since in the paper we are specifically interested in the case where \vector{v} and the transformed
vector are real. We then show how the left and right Householder transformations are related.

In the appendix, we present verification of several important properties of the Householder matrix
in the quaternion case. Since these properties generalize from the complex case, their proofs are
not vital here, but they are included in the paper for completeness.

The authors of~\cite{Bunse-Gerstner:1989} used quaternion Householder transformations in their \qr
algorithm and gave an algorithm for computing a Householder vector, based on that of Morrison.
However, they did not discuss the difference between left and right Householder transformations.


\subsection{Left Householder vector and matrix}
\label{lefthouseholder}
In this case the Householder matrix \matrix{H} is applied to the column vector \vector{a} on the
left, and the scalar $z$ is multiplied on the left of the matrix
$\matrix{I}-\vector{u}\hermitian{\vector{u}}$.Thus
$\matrix{H}=z\left(\matrix{I}-\vector{u}\hermitian{\vector{u}}\right)$.

\begin{theorem}
\label{theoremlefthouseholder}
Given an arbitrary quaternion vector $\vector{a}\in\H^r$, and a real vector $\vector{v}\in\R^r$ with
unit norm; there exist a quaternion vector $\vector{u}\in\H^r$ with norm
$\norm{\vector{u}}=\sqrt{2}$ and a unit quaternion scalar $z$, such that $\matrix{H}\vector{a} =
z\left(\matrix{I}-\vector{u}\hermitian{\vector{u}}\right)\vector{a} = \norm{\vector{a}}\vector{v}$.
\end{theorem}

\begin{proof}
The proof is by construction, and is based on that of Morrison~\cite{DOI:10.1145/321021.321030} for
the complex case, with the necessary adjustments to allow for non-commutative multiplication. We use
Morrison's symbols, and therefore we write $\zeta = \inverse{z}$ and $\alpha=\norm{\vector{a}}$.
Therefore, we have:
\begin{align}
\left(\matrix{I}-\vector{u}\hermitian{\vector{u}}\right)\vector{a} &= \zeta\vector{v}\alpha
\intertext{Multiplying out the parentheses on the left, we obtain:}
\label{luequals}
\vector{a}-\vector{u}\hermitian{\vector{u}}\vector{a} &= \zeta\vector{v}\alpha
\intertext{and then multiplying from the left by \hermitian{\vector{a}}, we obtain:}
\hermitian{\vector{a}}\vector{a}-
\hermitian{\vector{a}}\vector{u}\hermitian{\vector{u}}\vector{a} &=
\hermitian{\vector{a}}\zeta\vector{v}\alpha
\end{align}
Now, $\hermitian{\vector{a}}\vector{a} = \alpha^2$ by Lemma \ref{vectornorm}, and
$\hermitian{\vector{a}}\vector{u}\hermitian{\vector{u}}\vector{a}$ is the product of a quaternion
$\hermitian{\vector{a}}\vector{u}$ with its conjugate. Denoting the quaternion\footnote{In fact, we
make the choice later that $\mu$ is real.} $\hermitian{\vector{u}}\vector{a}$ by $\mu$, we then
have:
\begin{equation}
\label{lsubstitutehere}
\alpha^2 - |\mu|^2 = \hermitian{\vector{a}}\zeta\vector{v}\alpha
\end{equation}
The left-hand side of this equation is real, and therefore so must be the right-hand side. Therefore
the product $\hermitian{\vector{a}}\zeta$ is real, since \vector{v} and $\alpha$ are both real. At
this point, Morrison was able to reorder the terms and solve for $\zeta$ by writing
$\hermitian{\vector{a}}\vector{v}$ as the product of a positive real (modulus) with a unit complex
number. We cannot do this because $\zeta$ does not commute with the other terms. However, we can
re-order terms by the quaternion conjugate rule ($\conjugate{q_1 q_2} =
\conjugate{q_2}\,\conjugate{q_1}$). Thus: $\hermitian{\vector{a}}\zeta\vector{v}=
\conjugate{\conjugate{\zeta}\transpose{\vector{a}}}\vector{v}$. Since the result here is real, we
may drop the overall conjugate and obtain\footnote{It is at this point that we are forced to impose
the restriction that \vector{v} be real, since otherwise we could not re-order the right-hand side
to put $\zeta$ on the left of $\transpose{\vector{a}}$.}:
$\hermitian{\vector{a}}\zeta\vector{v}=\conjugate{\zeta}\transpose{\vector{a}}\vector{v}$. Following
Morrison, we now represent $\transpose{\vector{a}}\vector{v}$ as the product of a real modulus with
a unit quaternion, and write it as $\transpose{\vector{a}}\vector{v} = r\omega$, where $r$ is the
modulus, and $\omega$ is a unit quaternion. Then, the fact that
$\conjugate{\zeta}\transpose{\vector{a}}\vector{v}$ is real implies that
$\conjugate{\zeta}=\pm\,\conjugate{\omega}$ since $\zeta$ has unit modulus. This implies that
$\zeta=\pm\,\omega$. Substituting this result into equation \ref{lsubstitutehere}, and remembering
that $\omega$ is a unit quaternion, we obtain:
\begin{equation}
\alpha^2 - |\mu|^2 = \pm\,\conjugate{\omega}r\omega\alpha = \pm\,\alpha r
\end{equation}
which gives on re-arrangement:
\begin{equation}
\label{musquared}
|\mu|^2 = \alpha^2 \mp\,\alpha r
\end{equation}
We now make an arbitrary choice, and choose $\mu$ to be real (we could have done this earlier),
hence we have: $\mu = \sqrt{\alpha(\alpha\mp\,r)}$. Of the two possibilities here, we follow
Morrison and choose to add the terms inside the parentheses (both are positive, and if we subtract
them we might get a small result in some cases). Therefore we have $\mu = \sqrt{\alpha(\alpha+r)}$
and the choice we have just made also implies that $\zeta=-\omega$. Finally, from equation
\ref{luequals} we can obtain \vector{u}:
\begin{align}
\vector{a}-\zeta\vector{v}\alpha &=\vector{u}\hermitian{\vector{u}}\vector{a}=\vector{u}\mu\\
\intertext{and remembering that we have chosen $\mu$ to be real:}
\vector{u} &= \frac{1}{\mu}\left(\vector{a}-\zeta\vector{v}\alpha\right)
\end{align}
The sequence of steps necessary to compute \vector{u} and $\zeta$ is
therefore as follows:
\begin{align}
\alpha &= \norm{\vector{a}}\\
\intertext{If $\alpha = 0$ choose $\vector{u}$ to be a zero vector and $\zeta=1$ (this causes
the Householder matrix \matrix{H} to be an identity matrix), otherwise continue:}
r     &= |\transpose{\vector{a}}\vector{v}|\\
\label{zetal}
\zeta &= \begin{cases} 1 &:\quad r = 0\\
                       - \frac{\transpose{\vector{a}}\vector{v}}{r} &:\quad r > 0
         \end{cases}\\
\mu   &= \sqrt{\alpha(\alpha + r)}\\
\vector{u} &= \frac{1}{\mu}\left(\vector{a}-\zeta\vector{v}\alpha\right)
\end{align}
It remains to show that $\norm{\vector{u}}=\sqrt{2}$. From Lemma~\ref{vectornorm},
$\norm{\vector{u}}^2 = \hermitian{\vector{u}}\vector{u}$, hence:
\begin{align}
\norm{\vector{u}}^2 &=
\frac{1}{\mu}\hermitian{\left(\vector{a}-\zeta\vector{v}\alpha\right)}
\frac{1}{\mu}           \left(\vector{a}-\zeta\vector{v}\alpha\right) =
\frac{1}{\mu^2}\left(\hermitian{\vector{a}}-\conjugate{\zeta}\transpose{\vector{v}}\alpha\right)
                          \left(\vector{a}-\zeta\vector{v}\alpha\right)\\
\intertext{Multiplying out the terms in parentheses, and re-ordering terms that commute:}
\norm{\vector{u}}^2 &= \frac{1}{\mu^2}\left(\hermitian{\vector{a}}\vector{a}
                             -\conjugate{\zeta}\transpose{\vector{v}}\vector{a}\alpha
                             -\hermitian{\vector{a}}\vector{v}\zeta\alpha
                             +\conjugate{\zeta}\zeta\transpose{\vector{v}}\vector{v}\alpha^2
                        \right)\\
\intertext{$\hermitian{\vector{a}}\vector{a}=\alpha^2$ by Lemma \ref{vectornorm},
$\conjugate{\zeta}\zeta = |\zeta|^2 = 1$ and $\transpose{\vector{v}}\vector{v}=\norm{\vector{v}}^2=1$, and
$\mu^2 = \alpha(\alpha+r)$, hence:}
\norm{\vector{u}}^2 &= \frac{\alpha^2 -\conjugate{\zeta}\transpose{\vector{v}}\vector{a}\alpha
                                      -\hermitian{\vector{a}}\vector{v}\zeta\alpha
                            +\alpha^2}{\alpha(\alpha+r)} =
                      \frac{2\alpha-\conjugate{\zeta}\transpose{\vector{v}}\vector{a}
                                   -\hermitian{\vector{a}}\vector{v}\zeta}{\alpha+r}\\
\intertext{Since \vector{v} is real, $\transpose{\vector{v}}\vector{a} = \transpose{\vector{a}}\vector{v}$,
and we have:}
\norm{\vector{u}}^2 &= \frac{2\alpha
                                  -\conjugate{\zeta}\transpose{\vector{a}}\vector{v}
                                  -\hermitian{\vector{a}}\vector{v}\zeta}{\alpha+r}\\
\intertext{Finally, from equation \ref{zetal}, $\transpose{\vector{a}}\vector{v} = -r\zeta$, and
therefore $\hermitian{\vector{a}}\vector{v} = -r\conjugate{\zeta}$ since \vector{v} and $r$ are real.
Therefore:}
\norm{\vector{u}}^2 &= \frac{2\alpha+r\conjugate{\zeta}\zeta+r\conjugate{\zeta}\zeta}{\alpha+r}
=\frac{2\alpha+2r}{\alpha+r}=2
\end{align}
as required.
\qquad\end{proof}


\subsection{Right Householder vector and matrix}
\label{righthouseholder}
In this case the Householder matrix \matrix{G} is applied to the row vector \transpose{\vector{a}}
on the right, and the scalar $z$ is multiplied on the right of the matrix
$\matrix{I}-\vector{u}\hermitian{\vector{u}}$. Thus
$\matrix{G}=\left(\matrix{I}-\vector{u}\hermitian{\vector{u}}\right)z$.

\begin{theorem}
Given an arbitrary quaternion vector $\transpose{\vector{a}}\in\H^r$, and a real vector
$\transpose{\vector{v}}\in\R^r$ with unit norm, there exist a quaternion vector
$\transpose{\vector{u}}\in\H^r$ with norm $\norm{\vector{u}}=\sqrt{2}$ and a unit quaternion scalar
$z$, such that $\transpose{\vector{a}}\matrix{G} =
\transpose{\vector{a}}\left(\matrix{I}-\vector{u}\hermitian{\vector{u}}\right)z =
\norm{\vector{a}}\transpose{\vector{v}}$.
\end{theorem}

\begin{proof}
This theorem follows from Theorem \ref{theoremlefthouseholder}:
\begin{align*}
\transpose{\vector{a}}\left(\matrix{I}-\vector{u}\hermitian{\vector{u}}\right)z &=
\norm{\vector{a}}\transpose{\vector{v}}
\intertext{Place $z$ on the left using the quaternion conjugate rule
           $\left(\conjugate{q_1 q_2} = \conjugate{q_2}\,\conjugate{q_1}\right)$:}
\conjugate{z}\,\conjugate{\transpose{\vector{a}}\left(\matrix{I}-\vector{u}\hermitian{\vector{u}}\right)} &=
\norm{\vector{a}}\transpose{\vector{v}}
\intertext{Transpose both sides:}
\conjugate{z}\,\hermitian{\transpose{\vector{a}}\left(\matrix{I}-\vector{u}\hermitian{\vector{u}}\right)} &=
\norm{\vector{a}}\vector{v}
\intertext{Using Lemma~\ref{matrixconjugatetranspose}, we place \vector{a} on the right:}
\conjugate{z}\,\hermitian{\left(\matrix{I}-\vector{u}\hermitian{\vector{u}}\right)}\conjugate{\vector{a}} &=
\norm{\vector{a}}\vector{v}
\intertext{Since $z$ is a scalar we may take it under the overall conjugate transpose:}
\hermitian{z\left(\matrix{I}-\vector{u}\hermitian{\vector{u}}\right)}\conjugate{\vector{a}} &=
\norm{\vector{a}}\vector{v}
\end{align*}
We recognise the left-hand side to be $\hermitian{\matrix{H}}\!\conjugate{\vector{a}}$ as in
Theorem~\ref{theoremlefthouseholder}. That is, we may compute a right Householder matrix
by taking the conjugate transpose of the left Householder matrix computed from a conjugated
(and transposed) vector \vector{a}.
\end{proof}

For the sake of completeness, we give without proof the algorithm required to compute a right
Householder transformation directly. The sequence of steps necessary to compute
\transpose{\vector{u}} and $\zeta$ is as follows:
\begin{align}
\alpha &= \norm{\vector{a}}\\
\intertext{If $\alpha = 0$ choose $\vector{u}$ to be a zero vector and $\zeta=1$, otherwise continue:}
r     &= |\transpose{\vector{v}}\vector{a}|\\
\zeta &= \begin{cases} 1 &:\quad r = 0\\
                       - \frac{\transpose{\vector{v}}\vector{a}}{r} &:\quad r > 0
         \end{cases}\\
\mu   &= \sqrt{\alpha(\alpha + r)}\\
\transpose{\vector{u}} &=
\frac{1}{\mu}\left(\alpha\transpose{\vector{v}}\conjugate{\zeta}-\hermitian{\vector{a}}\right)
\end{align}


\section{Existence of the bidiagonalization of a quaternion matrix}
\label{existence}
Golub and Kahan in their 1965 paper~\cite{GolubKahan:1965} demonstrated the existence of the
bidiagonalization of a real matrix and showed how to decompose an arbitrary matrix into the product
of two orthogonal matrices and a bidiagonal matrix, either using the method now known as the
Golub-Kahan method, or by using a sequence of Householder transformations. We show here, using
quaternion Householder transformations, that the same decomposition is valid for an
arbitrary quaternion matrix and, more significantly, that the bidiagonal matrix may be real. This
latter property follows from the behaviour of the quaternion Householder transformations in
§\ref{quaternionhouseholder} when the vector \vector{v} is real.

\begin{theorem}
\label{existencetheorem}
Given an arbitrary quaternion matrix $\matrix{A}\in\H^{r×c}$ with $r$ rows and $c$ columns, there
exists a pair of unitary quaternion matrices $\matrix{L}\in\H^{r×r}$ and $\matrix{R}\in\H^{c×c}$,
and a \emph{real} bidiagonal matrix $\matrix{B}\in\R^{r×c}$ such that
$\matrix{L}\matrix{A}\matrix{R} = \matrix{B}$. \matrix{B} is upper or lower bidiagonal dependent on
the relative magnitudes of $r$ and $c$. Specifically, if $c>r$ ($r>c$), \matrix{B} is non-zero only
on its diagonal and super-diagonal (sub-diagonal respectively). If $r=c$, \matrix{B} may be non-zero
on either its sub- or super-diagonal arbitrarily.
\end{theorem}
\begin{proof}
The proof is by induction and construction. We show that the first column of \matrix{A} may
be transformed so that the first element is the only non-zero element and is real. We then show that
the remainder of the first row may be similarly transformed by transposing and conjugating the
remaining columns and applying the same process to the first column of the result. The process
continues until the matrix to be transformed has one column only. Transforming this column completes
the process of bidiagonalization.

Assume for the moment that $c>r$ and therefore that \matrix{B} is to be upper bidiagonal.
Figure~\ref{colrows} shows how the upper bidiagonal matrix \matrix{B} may be constructed recursively
from columns and rows with only one non-zero element.
\begin{figure}[tbp]
%
%
\newcommand{\nz}{\ensuremath ×}
\newcounter{x}
\newcounter{y}
\newcounter{s}
\newcounter{t}
\newcounter{k}
\setlength{\unitlength}{1.4em}
\newcommand{\bidiagonal}[2]{\setcounter{x}{#2}%
                            \setcounter{y}{#1}%
                            \begin{picture}(\thex,\they)
                            \setcounter{t}{\they-1}
                            \setcounter{s}{\thex-2}
                            \setcounter{k}{\they-3}
                            \multiput(0,2)(0,1){\thek}{\makebox(1,1){0}}
                            \put(0,1){\makebox(1,1){\vdots}}
                            \put(0,0){\makebox(1,1){0}}
                            \put(0,\thet){\makebox(1,1){\nz}}
                            \setcounter{k}{4*\they}
                            \ifthenelse{\thex>1}{\multiput(1,0)(0,0.25){\thek}{\line(0,1){0.125}}
                                                 \put(1,\thet){\makebox(1,1){\nz}}
                                                 \multiput(2,\thet)(1,0){\thes}{\makebox(1,1){0}}
                                                 \setcounter{k}{4*(\thex-1)}
                                                 \multiput(1,\thet)(0.25,0){\thek}{\line(1,0){0.125}}
                                                 \addtocounter{x}{-1}
                                                 \addtocounter{y}{-1}
                                                 \put(1,0){\bidiagonal{\they}{\thex}}
                                                }{\relax}
                            \end{picture}}
\begin{center}
\setcounter{x}{8}
\setcounter{y}{10}
\begin{picture}(\thex,\they)
\put(0,0){\line(0,1){\they}}\put(0,0){\line(1,0){0.25}}\put(0,\they){\line(1,0){0.25}}
\put(\thex,0){\line(0,1){\they}}\put(\thex,0){\line(-1,0){0.25}}\put(\thex,\they){\line(-1,0){0.25}}
\put(0,0){\bidiagonal{\they}{\thex}}
\end{picture}
\end{center}
\caption{\label{colrows}Construction of an upper bidiagonal matrix from columns
                        and rows each with one non-zero leading element.
                        (Non-zero elements are marked by \nz.)}
\end{figure}
We may transform \matrix{A} using a left Householder transformation as defined in
Theorem~\ref{theoremlefthouseholder} so that all elements of the first column become zero, except
for the first element which becomes real (that is, it has zero vector part). The required left
Householder matrix \matrix{H} is computed by choosing $\vector{a}$ equal to the first column of
\matrix{A}, and $\vector{v} = \transpose{(1,0,0,0,\ldots,0)}$ of length $r$, the same length as
\vector{a}. The left unitary matrix \matrix{L} is then equal to \matrix{H}. Multiplying \matrix{A}
on the left by \matrix{L} transforms the first column (and modifies all the subsequent columns).
Since \vector{v} is real and has only one non-zero element, the transformed first column of
\matrix{A} will also be real with one non-zero element as shown in Figure~\ref{colrows}. The right
unitary matrix \matrix{R} is an identity matrix.

If \matrix{A} has only one column, it is now in (degenerate) real bidiagonal form as is evident from
the first column of Figure~\ref{colrows}.

If \matrix{A} has more than one column, let $\matrix{A}^\prime$ be the matrix with $r$ rows and
$c-1$ columns obtained by deleting the first column of \matrix{A}. $\matrix{A}^\prime$ may be
bidiagonalized by applying the above process to $\hermitian{\matrix{A}^\prime}$ by
Lemma~\ref{matrixconjugatetranspose}. The result will be:
$\hermitian{\matrix{R}^\prime}\hermitian{\matrix{A}^\prime}\hermitian{\matrix{L}^\prime} =
\hermitian{\matrix{B}^\prime}$ where $\matrix{L}^\prime$ will be the same size as \matrix{L}, and
$\matrix{R}^\prime$ will be one row and column smaller than \matrix{R}.
Replacing $\matrix{A}^\prime$ by $\matrix{B}^\prime$ and multiplying \matrix{L} on the left
by $\matrix{L}^\prime$ and replacing the submatrix $\matrix{R}(2\to c,2\to c)$ by
$\matrix{R}^\prime$ completes the proof.

In the case when $r>c$, we have from Lemma~\ref{matrixconjugatetranspose}:
\[
\hermitian{\matrix{R}}\hermitian{\matrix{A}}\hermitian{\matrix{L}} = \hermitian{\matrix{B}}
\]
and thus \matrix{L}, \matrix{R} and \matrix{B} exist because we may form \hermitian{\matrix{A}}.
Since \matrix{B} is real, $\hermitian{\matrix{B}} = \transpose{\matrix{B}}$

In the case where \matrix{A} is square, and $r=c$, we may choose \matrix{B} to be upper or lower
bidiagonal arbitrarily. The simplest choice is upper bidiagonal.
\qquad\end{proof}

A dual proof is possible using right Householder transformations applied to rows of \matrix{A}. In
this case the result would be lower bidiagonal, and the upper bidiagonal case would be handled by
applying the process to $\hermitian{\matrix{A}}$.

We note, incidentally, that since real and complex numbers may be treated as subsets of the
quaternions, the above proof is also valid for real and complex matrices (of course, the result is
well-known in these cases).


\section{Bidiagonalization using Householder transformations}
\label{bidiag}
We give here a concise recursive algorithm to bidiagonalize an arbitrary quaternion matrix based on
left transformations only in which rows are transformed by conjugating and transposing the matrix.
Algorithm \ref{bidiagalg} is a reference algorithm based on explicit Householder matrices. It is
easily implemented (we have done so in Matlab, as part of our quaternion library \cite{qtfm}),
but for large matrices, or matrices in which one dimension is large, the explicit Householder matrix
may be too large to be computed explicitly.

\renewcommand{\lstlistingname}{Algorithm}
\begin{lstlisting}[float=p,print=true,frame=lines,caption={Quaternion bidiagonalization},label=bidiagalg]
Input  : $\matrix{A}\in\H^{r×c}$   
Output : $\matrix{L}\in\H^{r×r}$, $\matrix{B}\in\R^{r×c}$, $\hermitian{\matrix{R}}\in\H^{c×c}$

if $c \le r$
  Bidiagonalize($\matrix{A}$, $\matrix{L}$, $\matrix{B}$, $\matrix{R}$)
else
  Bidiagonalize($\hermitian{\matrix{A}}$, $\matrix{R}$, $\matrix{B}$, $\matrix{L}$); $\matrix{B} = \transpose{\matrix{B}}$
end

procedure Bidiagonalize(Input  : $\matrix{A}\in\H^{r×c}$   
                        Output : $\matrix{L}\in\H^{r×r}$, $\matrix{B}\in\R^{r×c}$, $\matrix{R}\in\H^{c×c}$)

!\textit{Transform the first column of \matrix{A}, where \vector{v} is the vector $\transpose{(1,0,0,0,\ldots)}$}!

$\matrix{L} = H(\matrix{A}(1\range r, 1), \vector{v})$,  $\matrix{A} = \matrix{L}\matrix{A}$             
$\matrix{R} = \matrix{I}^{c×c}$

!\textit{If there is more than one column, transform the remainder of the}!
!\textit{first row by applying the procedure recursively to the conjugate}!
!\textit{transpose of the remaining columns:}!

if $c>1$
  Bidiagonalize($\hermitian{\matrix{A}(1\range r,2\range c)}$, $\matrix{R}^\prime$, $\matrix{T}$, $\matrix{L}^\prime$)
  $\matrix{L} = \matrix{L}^\prime\,\matrix{L}$ 
  $\matrix{A}(1\range r,2\range c) = \transpose{\matrix{T}}$
  $\matrix{R}(2\range c,2\range c) = \matrix{R}^\prime$
end if
$\matrix{B} =\,\text{Real part of}\,\matrix{A}$  !\textit{(The vector part of \matrix{A} will be zero.)}!

function H($\vector{a}\in\H^m$, $\vector{v}\in\R^m$)
  return !$\matrix{H}\in\H^{m×m}$ as defined in §\ref{lefthouseholder} using $\vector{a}$ and $\vector{v}$.!
end
\end{lstlisting}


\begin{lstlisting}[float=p,print=true,frame=lines,caption={Quaternion singular value decomposition},label=qsvdalg]
Input  : $\matrix{A}\in\H^{r×c}$   
Output : $\matrix{U}\in\H^{r×r}$, $\matrix{\Sigma}\in\R^{r×c}$, $\matrix{V}\in\H^{c×c}$

!Bidiagonalize \matrix{A} using Algorithm \ref{bidiagalg}, to obtain \matrix{L}, \matrix{B} and \hermitian{\matrix{R}}, such that!
$\matrix{L}\matrix{B}\hermitian{\matrix{R}}=\matrix{A}$

!Compute the SVD of \matrix{B}, to obtain $\matrix{W}\in\R^{r×r}$, \matrix{\Sigma}, $\matrix{X}\in\R^{c×c}$, such that!
$\matrix{B} = \matrix{W}\matrix{\Sigma}\transpose{\matrix{X}}$

$\matrix{U} = \hermitian{\matrix{L}}\matrix{W} \qquad \hermitian{\matrix{V}} = \transpose{\matrix{X}}\hermitian{\matrix{R}}$
\end{lstlisting}


\section{Quaternion Singular Value Decomposition}
\label{qsvd}
The existence of the bidiagonalization of a quaternion matrix is not in itself sufficient to prove
that computing the singular values of the bidiagonal matrix yields the singular values of the original
quaternion matrix. We therefore prove that the singular values of the bidiagonal real matrix are the
same as the singular values of the original quaternion matrix.
\begin{theorem}
Given an arbitrary quaternion matrix $\matrix{A}\in\H^{r×c}$, and a real bidiagonal matrix
$\matrix{B}\in\R^{r×c}$ as defined in Theorem~\ref{existencetheorem}, the singular values of
\matrix{A} are the same as the singular values of \matrix{B}.
\end{theorem}
\begin{proof}
From Theorem \ref{existencetheorem} there exist unitary quaternion matrices \matrix{L} and \matrix{R} that
will transform \matrix{A} to \matrix{B}, that is $\matrix{L}\matrix{A}\matrix{R} = \matrix{B}$ and
since \matrix{L} and \matrix{R} are unitary, $\matrix{A} =
\hermitian{\matrix{L}}\matrix{B}\hermitian{\matrix{R}}$. The singular value decomposition of
$\matrix{B} = \matrix{U}\matrix{\Sigma}\transpose{\matrix{V}}$ where \matrix{U} and \matrix{V} are
orthogonal, hence $\matrix{A} =
\hermitian{\matrix{L}}\matrix{U}\matrix{\Sigma}\transpose{\matrix{V}}\hermitian{\matrix{R}}$. From
the uniqueness of the singular values, and from the fact that $\hermitian{\matrix{L}}\matrix{U}$ is
unitary and $\transpose{\matrix{V}}\hermitian{\matrix{R}}$ is unitary, it follows that
\matrix{\Sigma} contains the singular values of the quaternion matrix \matrix{A}.
\qquad\end{proof}

Algorithm \ref{qsvdalg} presents formally the steps needed to compute the quaternion singular value
decomposition using the bidiagonalization presented in §\ref{bidiag}. The \svd of the real matrix
resulting from the algorithm presented in §\ref{bidiag} may be computed using any existing real
\svd algorithm or routine, for example, the \svd function in Matlab, or a suitable \lapack
subroutine (note that there are specialised routines in \lapack for computing the \svd of real
bidiagonal matrices, and one of these should be used in preference to an \svd routine for general
real matrices to exploit the structure of the bidiagonal matrix).


\section{Conclusions}
We have presented a practical and efficient method for computing the singular value decomposition of
an arbitrary quaternion matrix \matrix{A}. This method bidiagonalizes \matrix{A} to a real
bidiagonal matrix \matrix{B} using quaternion Householder transformations and a simple recursive
algorithm based on that of Golub and Kahan~\cite{GolubKahan:1965}. The \svd of \matrix{B}, computed
with any available \svd algorithm for real matrices, followed by the multiplication of each of two
quaternion matrices by a real matrix, yields the \svd of \matrix{A}.

The few previous methods for computing the quaternion \svd were based on complex adjoint matrices, or
direct computation of the quaternion \svd itself. The former suffers from the serious disadvantage
that the computation is not accurate, and the singular values and vectors are computed twice over.
The latter suffers from the disadvantage that it is complex to code, and requires a deep
understanding of the \svd. To achieve the same computational efficiency and accuracy as the best
available real and complex \svd routines would require considerable effort in coding and testing
which the method presented in this paper neatly circumvents. In its place, it is only necessary
to code an efficient version of algorithms~\ref{bidiagalg} and \ref{qsvdalg}.


\appendix
\section{Appendix}
\label{elementary}
We gather here some elementary results needed in the paper. Most of these results would be trivial
in the complex case, but because of non-commutative quaternion multiplication we must be careful.

Elementary operations on quaternion matrices have been well-studied and gathered together in Zhang's
1997 paper~\cite{Zhang:97}. These include the concepts of unitary and Hermitian matrices which
generalize from complex matrices to quaternion matrices without difficulty. The product of unitary
quaternion matrices is itself unitary.

\begin{lemma}
\label{matrixconjugatetranspose}
For any two quaternion matrices of compatible dimensions:
$\hermitian{\matrix{A}\matrix{B}}=\hermitian{\matrix{B}}\hermitian{\matrix{A}}$.
\end{lemma}
\begin{proof}
This result was given by Zhang~\cite[Theorem 4.1 (2)]{Zhang:97}.
\end{proof}

\begin{lemma}
\label{vectornorm}
Given an arbitrary quaternion vector \vector{u}, its norm is given by
$\sqrt{\hermitian{\vector{u}}\vector{u}}$ or by
$\sqrt{\transpose{\vector{u}}\conjugate{\vector{u}}}$.
\end{lemma}
\begin{proof}
The norm of a vector $\norm{\vector{u}} = \sqrt{\sum_i |u_i|^2}$.
The inner product of \vector{u} with its transpose conjugate yields the sum of the elements of
\vector{u} with their conjugates. Since a quaternion multiplied by its conjugate gives the square
of its modulus, the inner product gives the sum of the moduli squared, and taking the square root
gives the norm of the vector.
\qquad\end{proof}

\begin{lemma}
\label{hermitianouterproduct}
Given an arbitrary quaternion vector \vector{u}, the matrix $\vector{u}\hermitian{\vector{u}}$ is
Hermitian.
\end{lemma}
\begin{proof}
Elements of the matrix on the diagonal consist of quaternions multiplied by their conjugates $u_i
\conjugate{u_i}$, and are therefore real. Elements not on the diagonal are quaternion conjugate
pairs, since they are of the form either $u_i \conjugate{u_j}$ or $\conjugate{u_i} u_j$.
\qquad\end{proof}

\begin{theorem}
A quaternion Householder matrix, \matrix{H}, of the form $\matrix{I}-\vector{u}\hermitian{\vector{u}}$,
is unitary if \vector{u} is a quaternion vector with $\norm{\vector{u}}=\sqrt{2}$.
\end{theorem}
\begin{proof}
It is sufficient to show that $\matrix{H}^2 = \matrix{I}$, since $\matrix{H} =
\hermitian{\matrix{H}}$ from Lemma \ref{hermitianouterproduct}, and subtracting an Hermitian matrix
from an identity matrix gives an Hermitian matrix (self-evident):
\begin{equation}
\matrix{H}^2 = \left(\matrix{I} - \vector{u}\hermitian{\vector{u}}\right)^2 =
\matrix{I}^2 - \vector{u}\hermitian{\vector{u}}\matrix{I}
             - \matrix{I}\vector{u}\hermitian{\vector{u}} +
\vector{u}\hermitian{\vector{u}}\vector{u}\hermitian{\vector{u}} =
\matrix{I} - 2\vector{u}\hermitian{\vector{u}} +
\vector{u}\left(\hermitian{\vector{u}}\vector{u}\right)\hermitian{\vector{u}}
\end{equation}
The term in the parentheses on the right is the square of the norm of \vector{u} by Lemma
\ref{vectornorm}, and since $\norm{\vector{u}}=\sqrt{2}$, we have:
\begin{equation}
\matrix{H}^2 = \matrix{I} - 2\vector{u}\hermitian{\vector{u}}
                          + 2\vector{u}\hermitian{\vector{u}} = \matrix{I}
\end{equation}
\qquad\end{proof}

\begin{theorem}
A unitary quaternion matrix \matrix{U} scaled by a quaternion $z$ with unit modulus is unitary
whether $z$ is multiplied on the left or right of \matrix{U}.
\end{theorem}
\begin{proof}
Consider first the case where $z$ is multiplied on the left of \matrix{U}:
\begin{align}
\label{scalarunitary}
\hermitian{(z\matrix{U})} &=
\hermitian{
\begin{pmatrix}
z u_{1,1} & z u_{1,2} & \hdots\\
z u_{2,1} & z u_{2,2} & \hdots\\
\vdots    & \vdots    & \ddots
\end{pmatrix}
} =
\begin{pmatrix}
\conjugate{z u_{1,1}} & \conjugate{z u_{2,1}} & \hdots\\
\conjugate{z u_{1,2}} & \conjugate{z u_{2,2}} & \hdots\\
\vdots                & \vdots                & \ddots
\end{pmatrix}\\
\intertext{and applying the quaternion conjugate rule
           $\left(\conjugate{q_1 q_2} = \conjugate{q_2}\,\conjugate{q_1}\right)$:}\notag
\hermitian{(z\matrix{U})} &=
\begin{pmatrix}
\conjugate{u_{1,1}}\,\conjugate{z} & \conjugate{u_{2,1}}\,\conjugate{z} & \hdots\\
\conjugate{u_{1,2}}\,\conjugate{z} & \conjugate{u_{2,2}}\,\conjugate{z} & \hdots\\
\vdots        & \vdots        & \ddots
\end{pmatrix} =
\hermitian{\matrix{U}}\conjugate{z}
\end{align}
Clearly, then:
\begin{align}
\label{unitary1}
\hermitian{\left(z\matrix{U}\right)}\left(z\matrix{U}\right)&=
\hermitian{\matrix{U}}\conjugate{z}z\matrix{U}=
\hermitian{\matrix{U}}|z|^2\matrix{U}=
\hermitian{\matrix{U}}\matrix{U} = \matrix{I}
\intertext{and}
\label{unitary2}
\left(z\matrix{U}\right)\hermitian{\left(z\matrix{U}\right)}&=
z\matrix{U}\hermitian{\matrix{U}}\conjugate{z} =
z\matrix{I}\conjugate{z} = |z|^2\matrix{I} = \matrix{I}
\end{align}

If $z$ is on the right of the matrix \matrix{U}, modification of equation \ref{scalarunitary} shows
easily that $\hermitian{\left(\matrix{U}z\right)} = \conjugate{z}\hermitian{\matrix{U}}$ and the
results of multiplying the scaled matrix by its transpose conjugate are as in equations
\ref{unitary1} and \ref{unitary2}.
\qquad\end{proof}



\end{document}